\documentclass{amsart}
\usepackage{graphicx} 
\usepackage{amsmath}
\usepackage{amssymb}

\newcommand{\R}{\mathbb{R}}

\title[Carnot-Carath\'{e}odory metrics in three dimensions]{Carnot-Carath\'{e}odory metrics associated to degenerate elliptic operators in three dimensions}
\author{Lyudmila Korobenko}
\address{Reed College\\
	Portland, Oregon, USA}
\email{korobenko@reed.edu}
\author{Florian Meister}
\address{Reed College\\
	Portland, Oregon, USA}
\email{florianjm@reed.edu}
\author{Olive Ross}
\address{Reed College\\
	Portland, Oregon, USA}
\email{odanr109@gmail.com}
\date{\today }

\begin{document}

\begin{abstract}
This note is a companion paper to \cite{Korobenko}. Here we generalize some of the results of \cite[Chapter 7]{Korobenko} to the case of a $3\times 3$ matrix function $A(x)\approx \mathrm{diag}\{1,f(x_1), g(x_1)\}$. More precisely, we make explicit calculations of the geodesics in the Carnot-Carath\'{e}odory space associated to $A$, and provide estimates on the Lebesgue measures of metric balls centered at the origin in that space.
\end{abstract}

	\maketitle

\section{Introduction}
One powerful approach to study regularity theory for degenerate elliptic equations is through properties of metric spaces associated to the operator, see for example \cite{Koskela}. The idea to use Carnot-Carath\'{e}odory metric spaces in place of $\R^n$ to perform the Moser iteration to establish H\"{o}lder continuity of weak solutions to degenereate elliptic equations goes back to Franchi and Lanconelli \cite{Franchi}, and many more results of this flavor have been obtained since then, see e.g. \cite{SaWh4, Korobenko, KoSa21, CrRo21, DiFaMoRo23}. Despite this, there are very few results that give any explicit description of these metric spaces, such as geodesic curves and metric balls. The works \cite{Mon17} and \cite{Korobenko} contain such calculations for metric spaces associated to specific classes of two-dimensional operators. In this paper we extend some of the results of \cite{Korobenko} to a class of three-dimensional operators.

The Moser iteration uses two key ingredients: a Sobolev inequality, valid for all compactly supported Lipschitz functions; and a Caccioppoli inequality, valid for weak solutions to a uniformly elliptic equation of the form $Lu\equiv \nabla\cdot A\nabla u=0$. As a result one obtains local boundedness of weak solutions. Using in addition a Poincar\'{e} inequality, one can obtain further regularity, namely, H\"{o}lder continuity of weak solutions. It turns out that this scheme can be adapted to operators that are not uniformly elliptic, such as degenerate elliptic operators. Fabes Kenig and Serapioni \cite{FKS82} identified four essential conditions necessary for performing the Moser iteration, namely, (1) the doubling condition; (2) uniqueness of the gradient; (3) $(q,p)$ Sobolev inequality with $q>p$; and (4) $(p,p)$ Poincar\'{e} inequality.
Franchi and Lanconelly \cite{Franchi} were the first to apply the Moser technique in the Carnot-Carath\'{e}odory setting. Numerous generalizations of their results followed \cite{GuLa, Saloff, CrRo21}, including a recent result by one of the authors and collaborators \cite{Korobenko} for infinitely degenerate operators. The difficulty in this latter case is that the metric balls are no longer doubling, and the classical Sobolev inequality no longer holds. A weaker Orlicz Sobolev inequality might be used instead, but its proof requires explicit description of metric balls, and estimates on their Lebesgue measures. 

In \cite{Korobenko} the authors consider an operator of the form $L=\nabla A\nabla$ with the $n\times n$ matrix $A$ satisfying
\[
A(x)\approx \mathrm{diag}\{1,1,\dots,1, f^{2}(x_1)\},
\]
which is a generalization of the two dimensional case of
\[
A(x)\approx \mathrm{diag}\{1, f^{2}(x_1)\}.
\]
For the Carnot-Carath\'{e}odory metric space associated to $A$ they give 1) an explicit description of geodesic curves; 2) the Jacobian determinant that allows integration over metric balls; 3) sharp estimates on Lebesgue measures of the metric balls. In this paper we generalize some of these results to the case of a $3\times 3$ matrix of the form
\[
A(x)\approx \mathrm{diag}\{1, f^{2}(x_1), g^{2}(x_1)\}.
\]
There seems to be some indication that the three dimensional case seems to be the key step to developing a theory for $n\times n$ dimensional matrices satisfying
\[
A(x)\approx \mathrm{diag}\{1, f_{1}^{2}(x_1),\dots,  f_{n-1}^{2}(x_1)\}.
\]
More precisely, in \cite{KS3} the authors prove a regularity result for an operator whose matrix has the above form. It turns out that this result roughly only depends on the maximum and minimum of the functions $f_1,\dots, f_{n-1}$. This suggests that the three dimensional matrices considered in this paper might be the model case that opens the door to the development of a complete $n$ dimensional theory.

The main results of this paper are the implicit equations for geodesic curves, and the upper and lower bound estimates for Lebesgue measures of metric balls. Unlike in the two dimensional case, the upper and lower bounds do not match in general, only under some extra assumptions on $f$ and $g$. So far we do not know if this mismatch is crucial for applications in regularity theory.

This paper is organized as follows. In Section \ref{sec:geodesics} we obtain implicit equations describing the geodesics starting at the origin. Section \ref{sec:measure} is devoted to upper and lower estimates on the Lebesgue measure of a subunit ball $B(0,R)$, centered at the origin. Finally, the appendix contains detailed calculations of the change of variable determinant $\partial(x,y,z)/\partial(r,\lambda, \mu)$.

\section{Geodesic equations}\label{sec:geodesics}
Following \cite[Chapter 7]{Korobenko} we have that the $A$-distance is given by
\begin{align*}
    dt^2=dx^2+\frac1{f^2(x)}dy^2+\frac1{g^2(x)}dz^2.
\end{align*}
This gives
\begin{align*}
    \left(\frac{dt}{dx}\right)^2=1+\frac1{f^2(x)}\left(\frac{dy}{dx}\right)^2+\frac1{g^2(x)}\left(\frac{dz}{dx}\right)^2,
\end{align*}
and the distance from the origin to a point with $x=x_0$ is thus
\begin{align*}
    t=\int_0^{x_0}\sqrt{1+\frac1{f^2(x)}\left(\frac{dy}{dx}\right)^2+\frac1{g^2(x)}\left(\frac{dz}{dx}\right)^2}dx.
\end{align*}
The goal then is to find $y(x)$ and $z(x)$ that minimize this distance, since the subunit metric is the infimum of admissible paths. This can be done using calculus of variations. Assume $y$ and $z$ are such that they minimize the distance Define $\Phi(\delta)$ to be the distance
\begin{align*}
    d\left((0,0,0,),(x_0,y(x_0)+\delta\alpha(x_0),z(x_0)+\delta\beta(x_0)\right),
\end{align*}
where $\alpha$ and $\beta$ are variations. By assumption, this will have a minimum at $\delta=0,$ and therefore $\frac{d\Phi(\delta)}{d\delta}=0.$ Using the known formula for the distance to a point $x_0$, 
\begin{align*}
    \Phi(\delta)=\int_0^{x_0}\sqrt{1+\frac1{f^2(x)}\left(\frac{dy}{dx}+\delta\frac{d\alpha}{dx}\right)^2+\frac1{g^2(x)}\left(\frac{dz}{dx}+\delta\frac{d\beta}{dx}\right)^2}dx
\end{align*}
and
\begin{align*}
    \frac{d\Phi(\delta)}{d\delta}=&\int_0^{x_0}\frac1{2\sqrt{1+\frac1{f^2(x)}\left(\frac{dy}{dx}+\delta\frac{d\alpha}{dx}\right)^2+\frac1{g^2(x)}\left(\frac{dz}{dx}+\delta\frac{d\beta}{dx}\right)^2}}\\
    &\times\left(\frac2{f^2(x)}\left(\frac{dy}{dx}+\delta\frac{d\alpha}{dx}\right)\frac{d\alpha}{dx}+\frac2{g^2(x)}\left(\frac{dz}{dx}+\delta\frac{d\beta}{dx}\right)\frac{d\beta}{dx}\right)dx
\end{align*}
This implies that
\begin{align*}
    \left.\frac{d\Phi(\delta)}{d\delta}\right|_{\delta=0}=&\int_0^{x_0}\frac1{2\sqrt{1+\frac1{f^2(x)}\left(\frac{dy}{dx}+\delta\frac{d\alpha}{dx}\right)^2+\frac1{g^2(x)}\left(\frac{dz}{dx}+\delta\frac{d\beta}{dx}\right)^2}}\\
    &\times\left(\frac{1}{f^2}\frac{dy}{dx}\frac{d\alpha}{dx}+\frac{1}{g^2}\frac{dz}{dx}\frac{d\beta}{dx}\right)dx,
\end{align*}
which gives
\begin{align*}
    0=\int_0^{x_0}\frac1{\sqrt{1+\frac1{f^2(x)}\left(\frac{dy}{dx}+\delta\frac{d\alpha}{dx}\right)^2+\frac1{g^2(x)}\left(\frac{dz}{dx}+\delta\frac{d\beta}{dx}\right)^2}}\left(\frac{1}{f^2}\frac{dy}{dx}\frac{d\alpha}{dx}+\frac{1}{g^2}\frac{dz}{dx}\frac{d\beta}{dx}\right)dx.
\end{align*}
Since this must be true for any variations $\alpha$ and $\beta,$ including $0,$ it follows that 
\begin{align*}
    0=\int_0^{x_0}\frac1{\sqrt{1+\frac1{f^2(x)}\left(\frac{dy}{dx}\right)^2+\frac1{g^2(x)}\left(\frac{dz}{dx}\right)^2}}\frac{1}{f^2}\frac{dy}{dx}\frac{d\alpha}{dx}dx
\end{align*}
and
\begin{align*}
    0=\int_0^{x_0}\frac1{\sqrt{1+\frac1{f^2(x)}\left(\frac{dy}{dx}\right)^2+\frac1{g^2(x)}\left(\frac{dz}{dx}\right)^2}}\frac{1}{g^2}\frac{dz}{dx}\frac{d\beta}{dx}dx.
\end{align*}
Integrating these formulas by parts,
\begin{align*}
    0=\left.\frac1{\sqrt{1+\frac1{f^2(x)}\left(\frac{dy}{dx}\right)^2+\frac1{g^2(x)}\left(\frac{dz}{dx}\right)^2}}\frac{1}{g^2}\frac{dz}{dx}\beta\right|_0^{x_0}\\-\int_0^{x_0}\frac{d}{dx}\left(\frac1{g^2}\frac{dz}{dx}\frac1{\sqrt{1+\frac1{f^2(x)}\left(\frac{dy}{dx}\right)^2+\frac1{g^2(x)}\left(\frac{dz}{dx}\right)^2}}\right)\beta dx.
\end{align*}
The left component will be $0$, since the variation must be $0$ at the beginning and end of the path, and so
\begin{align*}
    0=\int_0^{x_0}\frac{d}{dx}\left(\frac1{g^2}\frac{dz}{dx}\frac1{\sqrt{1+\frac1{f^2(x)}\left(\frac{dy}{dx}\right)^2+\frac1{g^2(x)}\left(\frac{dz}{dx}\right)^2}}\right)\beta dx.
\end{align*}
Because this is true for any $\beta$,
\begin{align*}
    0=\frac{d}{dx}\left(\frac1{g^2}\frac{dz}{dx}\frac1{\sqrt{1+\frac1{f^2(x)}\left(\frac{dy}{dx}\right)^2+\frac1{g^2(x)}\left(\frac{dz}{dx}\right)^2}}\right).
\end{align*}
Therefore, $\frac1{g^2}\frac{dz}{dx}\frac1{\sqrt{1+\frac1{f^2(x)}\left(\frac{dy}{dx}\right)^2+\frac1{g^2(x)}\left(\frac{dz}{dx}\right)^2}}$ is conserved along the geodesic. Call this quantity $\mu.$ Likewise, via an identical argument for the $y$ component, $\frac1{f^2}\frac{dy}{dx}\frac1{\sqrt{1+\frac1{f^2(x)}\left(\frac{dy}{dx}\right)^2+\frac1{g^2(x)}\left(\frac{dz}{dx}\right)^2}}$ is also conserved. Call this quantity $\lambda.$
The definitions of the constants $\lambda$ and $\mu$ can be rearranged to show that
\begin{align*}
    \frac{dy}{dx}=\pm\lambda f^2\sqrt{\frac{1+\frac1{g^2}\left(\frac{dz}{dx}\right)^2}{1-\lambda^2f^2}}
\end{align*}
and
\begin{align*}
    \frac{dz}{dx}=\pm\mu g^2\sqrt{\frac{1+\frac1{f^2}\left(\frac{dy}{dx}\right)^2}{1-\mu^2g^2}}.
\end{align*}
This can be used to show that
\begin{align*}
    \frac{dy}{dx}=\pm\lambda f^2\sqrt{\frac{1}{1-\lambda^2f^2-\mu^2g^2}}
\end{align*}
and
\begin{align*}
    \frac{dz}{dx}=\pm\mu g^2\sqrt{\frac{1}{1-\lambda^2f^2-\mu^2g^2}}.
\end{align*}
Finally, integrating these with respect to $x$, we obtain
\begin{align*}
    y=\int_0^{x}\lambda f^2(t)\sqrt{\frac{1}{1-\lambda^2f^2(t)-\mu^2g^2(t)}}dt
\end{align*}
and
\begin{align*}
    z=\int_0^{x}\mu g^2(t)\sqrt{\frac{1}{1-\lambda^2f^2(t)-\mu^2g^2(t)}}dt.
\end{align*}
The figure below shows a geodesic starting at the origin. The point on the curve has distance $1$ from the origin.
\begin{figure}[htb!]
	\includegraphics[scale=0.9]{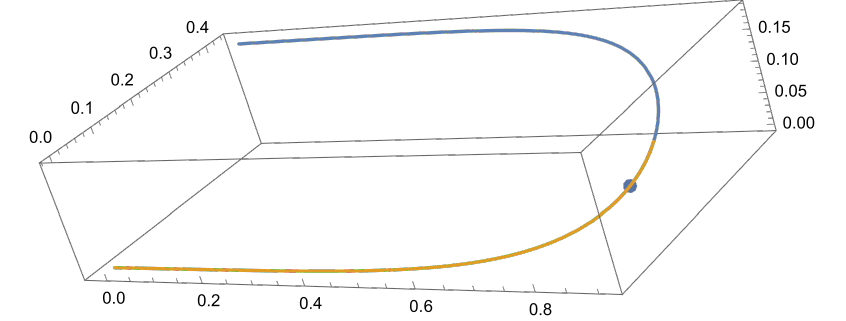}
	\caption{Geodesic starting from the origin}
\end{figure}

\section{Measure estimates}\label{sec:measure}
In this section we derive upper and lower bounds on Lebesgue measure of subunit balls centered at the origin.
\subsection{Upper bound}
For an upper bound estimate we first show that the projection of the three dimensional ball $B(0,R)$ on the $xy$-plane is contained in the two dimensional ball $B_{2D}(0,R)$, subunit with respect to the metric given by $\{1,f^2(x_1)\}$.
Let $L$ be a geodesic of length $R$ with the parametrization $\gamma:[0,R]\to \mathbb{R}^3$. Let $P$ be the projection of $L$ onto the $xy$-plane given by $\gamma^*:[0,R]\to \mathbb{R}^3$ where $\gamma^*(t)=(\gamma_{1}(t),\gamma_{2}(t),0)$.
Since $L$ is a geodesic, we know $\gamma$ satisfies the admissibility condition: for all $\xi\in\mathbb{R}^3$ we have
$$
(\xi_{1}\gamma'_{1}(t)+\xi_{2}\gamma'_{2}(t)+\xi_{3}\gamma'_{3}(t))^2\leq \xi_{1}^2+\xi_{2}^2f(\gamma_{1}(t))^2+\xi^3g(\gamma_{1}(t))^2
$$
This is true for all $\xi\in\mathbb{R}^3$ so therefore it is true for all $\xi$ of the form $(\xi_{1},\xi_{2},0)$ where $\xi_{1},\xi_{2}\in\mathbb{R}$. This means the following holds for all $\xi_{1},\xi_{2}\in\mathbb{R}$.
$$
(\xi_{1}\gamma'_{1}(t)+\xi_{2}\gamma'_{2}(t))^2\leq \xi_{1}^2+\xi_{2}^2f(\gamma_{1}(t))^2
$$
From the definition of $\gamma^*$ we obtain
\begin{align*}
(\xi_{1}(\gamma^*_{1})'(t)+\xi_{2}(\gamma^*_{2})'(t))^2 & \leq \xi_{1}^2+\xi_{2}^2f(\gamma^*_{1}(t))^2 \\
& \leq \xi_{1}^2+\xi_{2}^2f(\gamma^*_{1}(t))^2+\xi_{3}^2g(\gamma^*_{1}(t))^2.
\end{align*}
Therefore $\gamma^*$ is subunit with length $R$. This means the distance from the origin to the end point of $\gamma^*$ is less than $R$, therefore the end point of the projection is contained within the 2D ball of radius $R$ in the $xy$-plane. 
\\
The same argument tells us that the endpoint of the projection is contained in the 2D ball of radius $R$ in the $xz$-plane. We conclude that the endpoint of $L$ is contained within the shape outlined above.

To estimate the measure $|B(0,R)|$ we note that from \cite[Conclusion 45]{Korobenko} we know that the area of a 2D ball centered at the origin is proportional to $\frac{f(R)}{|F'(R)|^2}$, and from Proposition 47 we know that 2D balls centered on the origin have a maximum height proportional to $\frac{f(R)}{|F'(R)|}$. By above, we know that the ball $B(0,R)$ is contained in the cylindrical shape outlined by the boundary of
the 2D ball in the $xy$-plane and extending in the positive and negative $z$-direction to the height of the 2D ball in the $xz$-plane. Since the area of the 2D ball in the $xy$-plane is proportional to  $\frac{f(R)}{|F'(R)|^2}$, and the maximal height of the 2D ball in the $xz$-plane is approximately $\frac{g(R)}{|G'(R)|}$, the volume of the resulting shape is given by $\frac{f(R)g(R)}{|F'(R)|^2|G'(R)|}$.
Interchanging the roles of $y$ and $z$ we obtain the upper bound of $\frac{f(R)g(R)}{|F'(R)||G'(R)|^2}$.
Combining the two estimates gives
\begin{equation}\label{upper}
|B(0,R)|\lesssim\min\left\{ \frac{f(R)g(R)}{|F'(R)|^2|G'(R)|}, \frac{f(R)g(R)}{|F'(R)||G'(R)|^2} \right\}.
\end{equation}

\subsection{Lower bound}

For a lower bound estimate we use subunit curves of maximum parameter $R$ that are not necessarily geodesics. Since each of these curves is not necessarily a geodesic, its length is less than or equal to its maximum parameter $R$, so it will be contained within the metric ball of radius $R$. This means that the shape containing only these subunit curves is contained within the metric ball, so the volume of the metric ball is bounded below by the volume of the shape.

These subunit curves will be made up of 3 straight lines, one parallel to each axis. The curves will be of the form
\[
\varphi(t)= \begin{cases}
(t,0,0) & 0\leq t\leq a \\
(a,(t-a)f(a),0) & a<t\leq b \\
(a,(b-a)f(a),(t-b)g(a)) & b<t\leq R
\end{cases}.
\]
Each portion of the curve is subunit so therefore the whole curve is. We also know that the curve is continuous. To find the volume of the shape we need to find the maximum $x$, $y$, and $z$ values for each curve. The maximum $x$ value ranges from $0$ to $R$. For a curve with maximum $x$ value of $x_0$, the $y$ value ranges from 0 to $(R-x_0)f(x_0)$, achieving this maximum value when $b=R$. For a curve with maximum $x$ value of $x_0$ and maximum $y$ value of $y_0$, we want to find the maximum $z$ value in terms of $x_0$ and $y_0$ instead of $a$ and $b$. Note that $y_0=(b-x_0)f(x_0)$ so $b=\frac{y_0}{f(x_0)}+x_0$. Therefore the $z$ value ranges from 0 to $\left( R-x_0-\frac{y_0}{f(x_0)} \right)g(x_0)$.

We can now find the volume of the shape by integrating, which establishes a lower bound for the volume of the metric ball. Let $M=\min\{\frac{1}{|F'(R)|},\ \frac{1}{|G'(R)|}\}$,

\begin{align*}
V = & \int_{0}^R \int_{0}^{(R-x)f(x)} \int_{0}^{\left( R-x-\frac{y}{f(x)} \right)g(x)} 1 \, dz   \, dy   \, dx  \\
= & \int_{0}^R \int_{0}^{(R-x)f(x)} \left( R-x-\frac{y}{f(x)} \right)g(x)   \, dy   \, dx \\  
= & \int_{0}^R \int_{0}^{(R-x)f(x)} (R-x)g(x)- y \frac{g(x)}{f(x)}   \, dy   \, dx \\
= & \int_{0}^R \left.\left((R-x)g(x)y-\frac{1}{2} \frac{g(x)}{f(x)}y^2\right)\right|_{y=0}^{(R-x)f(x)} \, dx  \\
= & \int_{0}^R (R-x)g(x)(R-x)f(x)-\frac{1}{2} \frac{g(x)}{f(x)}(R-x)^2f(x)^2 \, dx  \\
= & \int_{0}^R (R-x)^2f(x)g(x) - \frac{1}{2}(R-x)^2f(x)g(x) \, dx  \\
= & \int_{0}^R \frac{1}{2}(R-x)^2f(x)g(x) \, dx  \\
\geq & \frac{1}{2}\int_{R-M}^R (R-x)^2f(x)g(x) \, dx  \\
\geq & \frac{1}{2}\int_{R-M}^R (R-x)^2f(R-M)g(R-M) \, dx \\
= & \frac{1}{2}f(R-M)\ g(R-M)\int_{R-M}^R (R-x)^2 \, dx  \\
\approx & \frac{1}{2}f(R) g(R) (x-R)^3|_{x=R-M}^R  \\
= & \frac{1}{2} f(R) g(R)[(R-R)^3-((R-M)-R)^3] \\
= & \frac{1}{2} f(R) g(R)[0+M]^3 \\
\approx & f(R)g(R)M^3.
\end{align*}
We therefore obtain
\begin{equation}\label{lower}
|B(0,R)|\gtrsim\frac{f(R)g(R)}{\max\{|F'(R)|,\ |G'(R)|\}^3}.
\end{equation}
Compared to the upper bound in (\ref{upper}), this quantity is strictly smaller in general. However, in a special case when $\frac{1}{|F'(R)|}\approx \frac{1}{|G'(R)|}$, the two bounds match. As an example of this one can consider $f(x)=e^{-1/x^2}$, $g(x)=e^{-2/x^2}$.

\section{Appendix}
We now calculate the determinant for the change of variables $(x,y,z)\to (r,\lambda,\mu)$.
In order to do this we must first find the derivatives of $x$, $y$, and $z$ with respect to $r$, $\lambda$, and $\mu$. This is primarily done using the chain rule for multiple variables, in order to utilize the derivatives we already know.

$\frac{dx}{dr}$:

\begin{align*}
r= & \int _{0}^x \frac{1}{\sqrt{ 1-\lambda^2f(u)^2-\mu^2g(u)^2 }} \, du  \\
\frac{dr}{dr}= & \frac{dr}{dx}\cdot \frac{dx}{dr} \\
1= & \left(\frac{1}{\sqrt{ 1-\lambda^2f(x)^2-\mu^2g(x)^2 }}\right)\cdot \frac{dx}{dr} \\
\frac{dx}{dr}= & \sqrt{ 1-\lambda^2f(x)^2-\mu^2g(x)^2 }
\end{align*}

$\frac{dx}{d\lambda}$:

\begin{align*}
r & =\int_{0}^x \frac{1}{\sqrt{ 1-\lambda^2f(u)^2-\mu^2g(u)^2 }} \, du \\
\frac{dr}{d\lambda}  & =\frac{dx}{d\lambda}\cdot\frac{dr}{dx}+\frac{dr}{d\lambda} \\
&=\frac{dx}{d\lambda}\left( \frac{1}{\sqrt{ 1-\lambda^2f(x)^2-\mu^2g(x)^2 }} \right) + \int_{0}^x \frac{ \partial  }{ \partial \lambda } \left[\frac{1}{\sqrt{ 1-\lambda^2f(u)^2-\mu^2g(u)^2 }}\right] \, du  \\
&=\frac{dx}{d\lambda}\left( \frac{1}{\sqrt{ 1-\lambda^2f(x)^2-\mu^2g(x)^2 }} \right) + \int_{0}^x \frac{\lambda f(u)^2}{ (1-\lambda^2f(u)^2-\mu^2g(u)^2)^{\ 3/2} } \, du \\
\frac{dx}{d\lambda} & = -\lambda \sqrt{ 1-\lambda^2f(x)^2-\mu^2g(x)^2 }\cdot \int_{0}^x \frac{f(u)^2}{ (1-\lambda^2f(u)^2-\mu^2g(u)^2)^{\ 3/2} } \, du
\end{align*}

$\frac{dx}{d\mu}$:

\begin{align*}
r = & \int_{0}^x \frac{1}{\sqrt{ 1-\lambda^2f(u)^2-\mu^2g(u)^2 }} \, du \\
\frac{dx}{d\mu}= & -\mu\sqrt{ 1-\lambda^2f(x)^2-\mu^2g(x)^2 }\cdot \int_{0}^x \frac{g(u)^2}{ (1-\lambda^2f(u)^2-\mu^2g(u)^2)^{\ 3/2} } \, du
\end{align*}

$\frac{dy}{dr}$:

\begin{align*}
y=&\int_{0}^x \frac{\lambda f(u)^2}{\sqrt{ 1-\lambda^2f(u)^2-\mu^2g(u)^2 }} \, du \\
\frac{dy}{dr}= & \frac{dy}{dx}\cdot \frac{dx}{dr} \\
= & \left(\frac{\lambda f(x)^2}{\sqrt{ 1-\lambda^2f(x)^2-\mu^2g(x)^2 }}\right)\cdot\sqrt{ 1-\lambda^2f(x)^2-\mu^2g(x)^2 } \\
= & \lambda f(x)^2
\end{align*}

$\frac{dy}{d\lambda}$:

\begin{align*}
y=&\int_{0}^x \frac{\lambda f(u)^2}{\sqrt{ 1-\lambda^2f(u)^2-\mu^2g(u)^2 }} \, du  \\
\frac{dy}{d\lambda} =&\frac{dy}{dx}\cdot \frac{dx}{d\lambda} +\frac{dy}{d\lambda} \\
=&\left( \frac{\lambda f(x)^2}{\sqrt{ 1-\lambda^2f(x)^2-\mu^2g(x)^2 }} \right)\cdot\left( -\lambda \sqrt{ 1-\lambda^2f(x)^2-\mu^2g(x)^2 }\int_{0}^x \frac{f(u)^2}{ (1-\lambda^2f(u)^2-\mu^2g(u)^2)^{\ 3/2} } \, du  \right)+ \\
&\int_{0}^x \frac{ \partial  }{ \partial \lambda } \left[\frac{\lambda f(u)^2}{\sqrt{ 1-\lambda^2f(u)^2-\mu^2g(u)^2 }}\right] \, du  \\
= & -\lambda^2f(x)^2\int_{0}^x \frac{f(u)^2}{ (1-\lambda^2f(u)^2-\mu^2g(u)^2)^{\ 3/2} } \, du + \int _{0}^x \frac{f(u)^2(1-\mu^2g(u)^2)}{(1-\lambda^2f(u)^2-\mu^2g(u)^2)^{\ 3/2}} \, du  \\
= & (1-\lambda^2f(x)^2)\int_{0}^x \frac{f(u)^2}{ (1-\lambda^2f(u)^2-\mu^2g(u)^2)^{\ 3/2} } \, du -\mu^2\int _{0}^x \frac{f(u)^2g(u)^2}{(1-\lambda^2f(u)^2-\mu^2g(u)^2)^{\ 3/2}} \, du
\end{align*}

$\frac{dy}{d\mu}$:

\begin{align*}
y=&\int_{0}^x \frac{\lambda f(u)^2}{\sqrt{ 1-\lambda^2f(u)^2-\mu^2g(u)^2 }} \, du \\
\frac{dy}{d\mu}= & \frac{dy}{dx}\cdot \frac{dx}{d\mu} + \frac{dy}{d\mu} \\
= & \left(\frac{\lambda f(x)^2}{\sqrt{ 1-\lambda^2f(x)^2-\mu^2g(x)^2 }} \right)\cdot\left(-\mu\sqrt{ 1-\lambda^2f(x)^2-\mu^2g(x)^2 }\cdot \int_{0}^x \frac{g(u)^2}{ (1-\lambda^2f(u)^2-\mu^2g(u)^2)^{\ 3/2} } \, du\right)+ \\
&\int_{0}^x \frac{ \partial  }{ \partial \mu } \left[\frac{\lambda f(u)^2}{\sqrt{ 1-\lambda^2f(u)^2-\mu^2g(u)^2 }}\right] \, du \\
= & -\lambda \mu f(x)^2\int_{0}^x \frac{g(u)^2}{ (1-\lambda^2f(u)^2-\mu^2g(u)^2)^{\ 3/2} } \, du + \lambda \mu \int _{0}^x \frac{f(u)^2g(u)^2}{(1-\lambda^2f(u)^2-\mu^2g(u)^2)^{\ 3/2}} \, du
\end{align*}

$\frac{dz}{dr}$:

\begin{align*}
z= & \int _{0}^x \frac{\mu g(u)^2}{\sqrt{ 1-\lambda^2f(u)^2-\mu^2g(u)^2 }} \, du \\
\frac{dz}{dr} = & \frac{dz}{dx}\cdot \frac{dx}{dr} \\
= & \left(\frac{\mu g(x)^2}{\sqrt{ 1-\lambda^2f(x)^2-\mu^2g(x)^2 }}\right)\cdot\sqrt{ 1-\lambda^2f(x)^2-\mu^2g(x)^2 } \\
= & \mu g(x)^2
\end{align*}

$\frac{dz}{d\lambda}$:

\begin{align*}
z= & \int _{0}^x \frac{\mu g(u)^2}{\sqrt{ 1-\lambda^2f(u)^2-\mu^2g(u)^2 }} \, du \\
\frac{dz}{d\lambda} = & \frac{dz}{dx}\cdot \frac{dx}{d\lambda} +\frac{dz}{d\lambda} \\
\dots& \\
= & -\lambda \mu g(x)^2\int_{0}^x \frac{f(u)^2}{ (1-\lambda^2f(u)^2-\mu^2g(u)^2)^{\ 3/2} } \, du + \lambda \mu \int _{0}^x \frac{f(u)^2g(u)^2}{(1-\lambda^2f(u)^2-\mu^2g(u)^2)^{\ 3/2}} \, du
\end{align*}

following the same argument as $\frac{dy}{d\mu}$.

$\frac{dz}{d\mu}$:

\begin{align*}
z= & \int _{0}^x \frac{\mu g(u)^2}{\sqrt{ 1-\lambda^2f(u)^2-\mu^2g(u)^2 }} \, du \\
\frac{dz}{d\mu} = & \frac{dz}{dx}\cdot \frac{dx}{d\mu} +\frac{dz}{d\mu} \\
\dots& \\
= & (1-\mu^2g(x)^2)\int_{0}^x \frac{g(u)^2}{ (1-\lambda^2f(u)^2-\mu^2g(u)^2)^{\ 3/2} } \, du -\lambda^2\int _{0}^x \frac{f(u)^2g(u)^2}{(1-\lambda^2f(u)^2-\mu^2g(u)^2)^{\ 3/2}} \, du
\end{align*}

following the same argument as $\frac{dy}{d\lambda}$.

Now we are ready to calculate the change of variable determinant, $\det \begin{bmatrix} \frac{dx}{dr} & \frac{dx}{d\lambda} & \frac{dx}{d\mu} \\ \frac{dy}{dr} & \frac{dy}{d\lambda} & \frac{dy}{d\mu} \\ \frac{dz}{dr} & \frac{dz}{d\lambda} & \frac{dz}{d\mu} \end{bmatrix}$.

\noindent The formula for this is $\frac{dx}{dr} \left(\frac{dy}{d\lambda}\cdot\frac{dz}{d\mu} -\frac{dy}{d\mu}\cdot \frac{dz}{d\lambda}\right) - \frac{dx}{d\lambda} \left( \frac{dy}{dr}\cdot\frac{dz}{d\mu} -\frac{dy}{d\mu}\cdot \frac{dz}{dr}\right) + \frac{dx}{d\mu} \left( \frac{dy}{dr}\cdot\frac{dz}{d\lambda} -\frac{dy}{d\lambda}\cdot \frac{dz}{dr}\right)$.

\noindent For legibility we set

\begin{align*}
A= & \int_{0}^x \frac{f(u)^2}{ (1-\lambda^2f(u)^2-\mu^2g(u)^2)^{\ 3/2} } \, du \\
B= & \int_{0}^x \frac{g(u)^2}{ (1-\lambda^2f(u)^2-\mu^2g(u)^2)^{\ 3/2} } \, du \\
C= & \int _{0}^x \frac{f(u)^2g(u)^2}{(1-\lambda^2f(u)^2-\mu^2g(u)^2)^{\ 3/2}} \, du
\end{align*}

Using the above substitutions we re-express our derivatives.

\begin{align*}
\frac{dx}{dr}= & \sqrt{ 1-\lambda^2f(x)^2-\mu^2g(x)^2 } \\
\frac{dx}{d\lambda}= & -\lambda \sqrt{ 1-\lambda^2f(x)^2-\mu^2g(x)^2 }\ A \\
\frac{dx}{d\mu}= & -\mu \sqrt{ 1-\lambda^2f(x)^2-\mu^2g(x)^2 }\ B \\
\\
\frac{dy}{dr}= & \lambda f(x)^2 \\
\frac{dy}{d\lambda}= & (1-\lambda^2f(x)^2)A -\mu^2C \\
\frac{dy}{d\mu}= & -\lambda \mu f(x)^2B + \lambda \mu C \\
\\
\frac{dz}{dr} = & \mu g(x)^2 \\
\frac{dz}{d\lambda}= & -\lambda \mu g(x)^2A + \lambda \mu C \\
\frac{dz}{d\mu}= & (1-\mu^2g(x)^2)B -\lambda^2C
\end{align*}

To help with legibility we will calculate the determinant one term at a time. The first term is $\frac{dx}{dr} \left(\frac{dy}{d\lambda}\cdot\frac{dz}{d\mu} -\frac{dy}{d\mu}\cdot \frac{dz}{d\lambda}\right)$. Plugging the above values into this expression gives the following.

\begin{align*}
& \sqrt{ 1-\lambda^2f(x)^2-\mu^2g(x)^2 }\ \bigg[((1-\lambda^2f(x)^2)A -\mu^2C)\cdot((1-\mu^2g(x)^2)B -\lambda^2C) - (-\lambda \mu f(x)^2B + \lambda \mu C)\cdot \\
&(-\lambda \mu g(x)^2A + \lambda \mu C)\bigg] \\
= & \sqrt{ 1-\lambda^2f(x)^2-\mu^2g(x)^2 }\ \bigg[(A-\lambda^2f(x)^2A -\mu^2C)\cdot(B-\mu^2g(x)^2B -\lambda^2C) -  \\
&(\lambda^2\mu^2f(x)^2g(x)^2AB - \lambda^2\mu^2g(x)^2AC - \lambda^2\mu^2f(x)^2BC + \lambda^2\mu^2C^2) \\ 
= & \sqrt{ 1-\lambda^2f(x)^2-\mu^2g(x)^2 }\ \bigg[(AB - \mu^2g(x)^2AB - \lambda^2AC - \lambda^2f(x)^2AB + \lambda^2\mu^2f(x)^2g(x)^2AB + \lambda^4f(x)^2AC  \\
&- \mu^2BC + \mu^4g(x)^2BC+\lambda^2\mu^2C^2) -  (\lambda^2\mu^2f(x)^2g(x)^2AB - \lambda^2\mu^2g(x)^2AC - \lambda^2\mu^2f(x)^2BC + \lambda^2\mu^2C^2) \\
= & \sqrt{ 1-\lambda^2f(x)^2-\mu^2g(x)^2 }\ \bigg[(AB - \mu^2g(x)^2AB - \lambda^2AC - \lambda^2f(x)^2AB + \lambda^4f(x)^2AC  \\
&- \mu^2BC + \mu^4g(x)^2BC) -  (- \lambda^2\mu^2g(x)^2AC - \lambda^2\mu^2f(x)^2BC))\bigg] \\
= & \sqrt{ 1-\lambda^2f(x)^2-\mu^2g(x)^2 }\ \bigg[AB - \mu^2g(x)^2AB - \lambda^2AC - \lambda^2f(x)^2AB + \lambda^4f(x)^2AC  \\
&- \mu^2BC + \mu^4g(x)^2BC + \lambda^2\mu^2g(x)^2AC + \lambda^2\mu^2f(x)^2BC\bigg]
\end{align*}

\noindent The second term to calculate is $-\frac{dx}{d\lambda} \left( \frac{dy}{dr}\cdot\frac{dz}{d\mu} -\frac{dy}{d\mu}\cdot \frac{dz}{dr}\right)$.

\begin{align*}
 & -(-\lambda \sqrt{ 1-\lambda^2f(x)^2-\mu^2g(x)^2 }\ A)\bigg[ \lambda f(x)^2\cdot((1-\mu^2g(x)^2)B -\lambda^2C) - (-\lambda \mu f(x)^2B + \lambda \mu C)\cdot\mu g(x)^2 \bigg] \\
= & \lambda \sqrt{ 1-\lambda^2f(x)^2-\mu^2g(x)^2 }\ A\bigg[ \lambda f(x)^2\cdot(B-\mu^2g(x)^2B -\lambda^2C) - (-\lambda \mu^2 f(x)^2g(x)^2B + \lambda \mu^2 g(x)^2C) \bigg] \\
= & \lambda \sqrt{ 1-\lambda^2f(x)^2-\mu^2g(x)^2 }\ A\bigg[ \lambda f(x)^2B-\lambda\mu^2f(x)^2g(x)^2B -\lambda^3f(x)^2C + \lambda \mu^2 f(x)^2g(x)^2B - \lambda \mu^2 g(x)^2C \bigg] \\
= & \lambda \sqrt{ 1-\lambda^2f(x)^2-\mu^2g(x)^2 }\ A\bigg[ \lambda f(x)^2B -\lambda^3f(x)^2C - \lambda \mu^2 g(x)^2C \bigg] \\
= & \sqrt{ 1-\lambda^2f(x)^2-\mu^2g(x)^2 }\ \bigg[ \lambda^2 f(x)^2AB -\lambda^4f(x)^2AC - \lambda^2 \mu^2 g(x)^2AC \bigg] 
\end{align*}

\noindent The final term is $\frac{dx}{d\mu} \left( \frac{dy}{dr}\cdot\frac{dz}{d\lambda} -\frac{dy}{d\lambda}\cdot \frac{dz}{dr}\right)$.

\begin{align*}
& -\mu\sqrt{ 1-\lambda^2f(x)^2-\mu^2g(x)^2 } B \bigg[ \lambda f(x)^2\cdot(-\lambda \mu g(x)^2A + \lambda \mu C) - ((1-\lambda^2f(x)^2)A -\mu^2C)\cdot\mu g(x)^2 \bigg]  \\
= & -\mu\sqrt{ 1-\lambda^2f(x)^2-\mu^2g(x)^2 } B\bigg[ -\lambda^2 \mu f(x)^2 g(x)^2A + \lambda^2 \mu f(x)^2 C - (A-\lambda^2f(x)^2A -\mu^2C)\cdot\mu g(x)^2 \bigg]  \\
= & -\mu\sqrt{ 1-\lambda^2f(x)^2-\mu^2g(x)^2 } B\bigg[ -\lambda^2 \mu f(x)^2 g(x)^2A + \lambda^2 \mu f(x)^2 C - (\mu g(x)^2A - \lambda^2 \mu f(x)^2 g(x)^2 A -\mu^3g(x)^2C) \bigg] \\
= & -\mu\sqrt{ 1-\lambda^2f(x)^2-\mu^2g(x)^2 } B\bigg[ \lambda^2 \mu f(x)^2 C - (\mu g(x)^2A - \mu^3g(x)^2C) \bigg] \\
= & -\mu\sqrt{ 1-\lambda^2f(x)^2-\mu^2g(x)^2 } B\bigg[ \lambda^2 \mu f(x)^2 C - \mu g(x)^2A + \mu^3g(x)^2C \bigg] \\
= & -\sqrt{ 1-\lambda^2f(x)^2-\mu^2g(x)^2 }\ \bigg[ \lambda^2 \mu^2 f(x)^2 BC - \mu^2 g(x)^2AB + \mu^4g(x)^2BC \bigg] 
\end{align*}

\noindent Combining together and simplifying we finally obtain

\begin{align*}
&\begin{bmatrix} \frac{dx}{dr} & \frac{dx}{d\lambda} & \frac{dx}{d\mu} \\ \frac{dy}{dr} & \frac{dy}{d\lambda} & \frac{dy}{d\mu} \\ \frac{dz}{dr} & \frac{dz}{d\lambda} & \frac{dz}{d\mu} \end{bmatrix}\\
&=  \sqrt{ 1-\lambda^2f(x)^2-\mu^2g(x)^2 }\ \bigg[AB - \mu^2g(x)^2AB - \lambda^2AC - \lambda^2f(x)^2AB + \lambda^4f(x)^2AC  \\
&- \mu^2BC + \mu^4g(x)^2BC + \lambda^2\mu^2g(x)^2AC + \lambda^2\mu^2f(x)^2BC\bigg] +  \\
&\sqrt{ 1-\lambda^2f(x)^2-\mu^2g(x)^2 }\ \bigg[ \lambda^2 f(x)^2AB -\lambda^4f(x)^2AC - \lambda^2 \mu^2 g(x)^2AC \bigg] -  \\
&\sqrt{ 1-\lambda^2f(x)^2-\mu^2g(x)^2 }\ \bigg[ \lambda^2 \mu^2 f(x)^2 BC - \mu^2 g(x)^2AB + \mu^4g(x)^2BC \bigg] \\
&=  \sqrt{ 1-\lambda^2f(x)^2-\mu^2g(x)^2 }\ \bigg[ AB - \mu^2g(x)^2AB - \lambda^2AC - \lambda^2f(x)^2AB + \lambda^4f(x)^2AC  \\
&- \mu^2BC + \mu^4g(x)^2BC + \lambda^2\mu^2g(x)^2AC + \lambda^2\mu^2f(x)^2BC +  \\
&\lambda^2 f(x)^2AB -\lambda^4f(x)^2AC - \lambda^2 \mu^2 g(x)^2AC -  \\
&(\lambda^2 \mu^2 f(x)^2 BC - \mu^2 g(x)^2AB + \mu^4g(x)^2BC) \bigg]  \\
&=  \sqrt{ 1-\lambda^2f(x)^2-\mu^2g(x)^2 }\ \bigg[ AB  - \lambda^2AC - \mu^2BC \bigg]
\end{align*}

\bibliographystyle{plain}
\bibliography{References}

\end{document}